\documentclass[a4paper,12pt,oneside,headsepline]{scrartcl}

\usepackage{etoolbox}
\usepackage{ifdraft}
\usepackage{iftex}
\usepackage{enumerate}

\ifluatex
\usepackage[utf8]{luainputenc}
\else
\usepackage[utf8]{inputenc}
\fi
\usepackage[T1]{fontenc}

\usepackage{lmodern}
\usepackage{fourier}

\usepackage{amsfonts}
\usepackage{mathrsfs} 
\usepackage{dsfont}

\usepackage{amssymb}
\usepackage{shuffle}

\usepackage{amsmath}
\usepackage{mathtools}
\usepackage[thmmarks, amsmath, amsthm]{ntheorem}
\usepackage{nccmath}    

\usepackage[draft=false]{scrlayer-scrpage}
\usepackage{needspace}

\usepackage[final,
            pdfborder={0 0 0}, colorlinks=true,
            linkcolor=BrickRed, citecolor=ForestGreen, urlcolor=RoyalBlue]{hyperref}

\usepackage[cmyk,dvipsnames]{xcolor}

\ifdraft{%
\usepackage[textsize=scriptsize,colorinlistoftodos]{todonotes}
\usepackage{lineno}
\usepackage{scrtime}
}{}

\usepackage{booktabs}
\usepackage[inline]{enumitem}
\usepackage{lettrine}
\usepackage{tikz-cd}
\usepackage{url}

\KOMAoptions{DIV=12}

\ifdraft{%
\linenumbers
\newcommand*\patchAmsMathEnvironmentForLineno[1]{%
  \expandafter\let\csname old#1\expandafter\endcsname\csname #1\endcsname
  \expandafter\let\csname oldend#1\expandafter\endcsname\csname end#1\endcsname
  \renewenvironment{#1}%
     {\linenomath\csname old#1\endcsname}%
     {\csname oldend#1\endcsname\endlinenomath}}%
\newcommand*\patchBothAmsMathEnvironmentsForLineno[1]{%
  \patchAmsMathEnvironmentForLineno{#1}%
  \patchAmsMathEnvironmentForLineno{#1*}}%
\AtBeginDocument{%
\patchBothAmsMathEnvironmentsForLineno{equation}%
\patchBothAmsMathEnvironmentsForLineno{align}%
\patchBothAmsMathEnvironmentsForLineno{flalign}%
\patchBothAmsMathEnvironmentsForLineno{alignat}%
\patchBothAmsMathEnvironmentsForLineno{gather}%
\patchBothAmsMathEnvironmentsForLineno{multline}%
}}

\clearmainofpairofpagestyles
\automark[section]{subsection}

\renewcommand{\subsectionmark}[1]{}

\lohead{{\small\normalfont \headertitle}}
\rohead{{\small \headerauthors}}

\RedeclareSectionCommand[%
  font=\Large\sffamily\bfseries,%
  beforeskip=1\baselineskip,%
  afterskip=0.5\baselineskip,%
  indent=0em,%
  afterindent=false%
  ]{section}

\RedeclareSectionCommands[%
  font=\normalfont\bfseries,%
  beforeskip=3pt,%
  afterskip=-1em%
  ]{subsection,subsubsection}

\RedeclareSectionCommands[%
  font=\normalfont\itshape,%
  beforeskip=1pt,%
  afterskip=-1em,%
  indent=0pt,%
  ]{paragraph}

\newenvironment{enumerateroman}{
\begin{enumerate}[label=(\roman*),%
  leftmargin=2.5em,itemindent=0pt,%
  labelindent=.5em,labelwidth=1.5em,labelsep=!,%
  nosep]
}{
\end{enumerate}
}

\newenvironment{enumeratearabic*}{
\begin{enumerate*}[label=(\arabic*)] 
}{
\end{enumerate*}
}

\newenvironment{enumerateroman*}{
\begin{enumerate*}[label=(\roman*)] 
}{
\end{enumerate*}
}

\numberwithin{equation}{section}

\theoremnumbering{arabic}
\newtheorem{theoremcounter}{theoremcounter}[section]
\theoremnumbering{Alph}

\theoremstyle{plain}

\newtheorem{conjecture}[theoremcounter]{Conjecture}
\newtheorem{corollary}[theoremcounter]{Corollary}
\newtheorem{lemma}[theoremcounter]{Lemma}

\newtheorem{proposition}[theoremcounter]{Proposition}
\newtheorem{theorem}[theoremcounter]{Theorem}

\newtheorem*{thm*}{Theorem}
\newtheorem*{cor*}{Corollary}

\theoremstyle{plain}

\theoremstyle{definition}

\newtheorem{definition}[theoremcounter]{Definition}
\newtheorem{example}[theoremcounter]{Example}

\theoremstyle{remark}

\newtheorem{remark}[theoremcounter]{Remark}

\theoremstyle{nonumberremark}

\newcommand{\tx}{\text}

\makeatletter
\newcommand{\writelabel}[1]{#1\def\@currentlabel{#1}}
\makeatother

\newcommand{\minwidthmathbox}[2]{%
  \mathmakebox[{\ifdim#1<\width\width\else#1\fi}]{#2}%
}

\newcommand{\tbf}{\bfseries}

\newcommand{\bbN}{\ensuremath{\mathbb{N}}}

\newcommand{\bbQ}{\ensuremath{\mathbb{Q}}}

\newcommand{\bbZ}{\ensuremath{\mathbb{Z}}}

\newcommand{\cH}{\ensuremath{\mathcal{H}}}

\newcommand{\rmH}{\ensuremath{\mathrm{H}}}

\newcommand{\mrelspace}[1]{\mathrel{\mspace{#1}}}

\let\rightarroworig\rightarrow
\renewcommand{\rightarrow}
  {\protect\relbar\mrelspace{-9.7mu}\rightarroworig}

\let\leftarroworig\leftarrow
\renewcommand{\leftarrow}
  {\protect\leftarroworig\mrelspace{-9.7mu}\relbar}

\renewcommand{\longrightarrow}
  {\protect\relbar\mrelspace{-3.2mu}\relbar\mrelspace{-9.5mu}\rightarroworig}

\makeatletter
\@ifundefined{xlongrightarrow}{%

}{%

}
\makeatother

\makeatletter
\@ifundefined{xlongleftarrow}{%

}{%

}
\makeatother

\newcommand{\Hom}{\operatorname{Hom}}
\newcommand{\id}{\mathrm{id}}

\newcommand{\NN}{\ensuremath{\mathbb{N}}}
\newcommand{\ZZ}{\ensuremath{\mathbb{Z}}}
\newcommand{\QQ}{\ensuremath{\mathbb{Q}}}

\DeclareMathOperator{\Ext}{Ext}

\DeclareMathOperator{\Mod}{Mod}

\DeclareMathOperator{\GrAlg}{GrAlg}
\DeclareMathOperator{\NAlg}{NUnAlg}
\DeclareMathOperator{\Sym}{Sym}

\DeclareMathOperator{\QS}{QS}

\newcommand{\headertitle}{{%
  $Z_\infty$-algebras: Zinbielification of $B_\infty$ and quasi-shuffle algebras
}}
\newcommand{\headerauthors}{%
  A.~Keilthy,
}
\title{%
  $Z_\infty$-algebras:\\ Zinbielification of $B_\infty$ and quasi-shuffle algebras
}
\author{%
Adam Keilthy%
}
\date{}

\begin{document}

\thispagestyle{scrplain}
\begingroup
\deffootnote[1em]{1.5em}{1em}{\thefootnotemark}
\maketitle
\endgroup


\begin{abstract}
\small
\noindent
{\tbf Abstract:}
Commutative $B_\infty$-algebras provide a natural generalisation of the shuffle and quasi-shuffle algebras, giving deformations of the shuffle algebra. We introduce a Zinbiel analogue in $Z_\infty$-algebras, giving deformations of the half-shuffle Zinbiel algebra. All commutative $B_\infty$-algebras have been shown to be isomorphic to the shuffle algebra, though this isomorphism is not canonical. We show that all $Z_\infty$-algebras are canonically isomorphic to the half-shuffle algebra, which gives a canonical isomorphism to the shuffle algebra for $B_\infty$-algebras arising from $Z_\infty$-algebras. Finally, we explore which $B_\infty$-algebras can be ``Zinbielified'' to a $Z_\infty$-algebra, working through a number of explicit examples in the case of quasi-shuffle algebras and uncovering a very explicit classification of simple Zinbielifications via semigroup cocycles.
\\[.3\baselineskip]
\noindent
\end{abstract}


\vspace{-1.5em}
\renewcommand{\contentsname}{}
\setcounter{tocdepth}{2}
\tableofcontents
\vspace{1.5em}


\addcontentsline{toc}{section}{Introduction}
\section*{Introduction}
A commutative $B_\infty$-algebra is a commutative Hopf algebra that is cofree as a conilpotent coalgebra. Explicitly, up to a choice of a basis the space of primitive elements, a commutative $B_\infty$-algebra is a commutative Hopf algebra on the vector space of non-commutative polynomials in variables indexed by this basis, such that the coproduct is given by deconcatenation. First introduced in the context of cochain complexes \cite{getzler94}, such algebras arise quite naturally in areas of algebra, combinatorics, and number theory. Classical examples include shuffle algebras \cite{Reutenauer1993} and quasi-shuffle algebras \cite{radford79,hoffman00,hoffmankentaro}, which underlie the algebraic structure of multiple zeta values \cite{Hoffman1997,Racinet2002}. Multiple zeta values also have ties to a more general $B_\infty$-algebra via the block shuffle algebra \cite{keilthy22, blockshuffle}. 

The shuffle algebra on an alphabet $X$ is referred to as the trivial $B_\infty$-algebra. It has underlying vector space $\QQ\langle X\rangle$, with product defined recursively by
\begin{align*}
    au \shuffle bv &:= a(u\shuffle bv) + b(au\shuffle v)\\
    w\shuffle 1 &= 1\shuffle w = w
\end{align*}
for all letters $a,b\in X$ and words $u,v,w\in X^\times$. A major result in the classification of commutative $B_\infty$-algebras due to Foissy and Patras \cite{binfty_foissy_patras,foissydeform} established that every commutative $B_\infty$-algebra on $\QQ\langle X\rangle$ is isomorphic to the shuffle algebra. Examples of such isomorphisms were known previously \cite{hoffman00,keilthy22}, though these were constructed on a case by case basis. 

However, these isomorphisms are not canonical, even if we restrict to isomorphisms that fix elements of $X$. Indeed, as $(\QQ\langle X\rangle,\shuffle)$ is isomorphic to a polynomial algebra in infinitely many variables, we can twist any such isomorphism by an automorphism of this polynomial algebra fixing only finitely many variables.

The half-shuffle Zinbiel algebra \cite{Loday1995} on $\QQ\langle X\rangle$ has product defined by
\begin{align*}
    y_iu \succ_\shuffle v &= y_i(u\shuffle v).\\
\end{align*}
In contrast to the shuffle algebra, the only Zinbiel automorphism of the half-shuffle algebra fixing $X$ is the identity. As such, if $(\QQ\langle X\rangle,\succ)$ is any other Zinbiel algebra, there is a \emph{unique} Zinbiel morphism fixing $X$ from the half-shuffle algebra. 

This rigidity provides the motivation for this article. We will defined a $Z_\infty$-algebra as a Zinbiel analogue of $B_\infty$-algebras, as cofree comodule over a cofree coalgebra. In this framework, we establish the following.
\begin{thm*}
    Given an vector space $V$, all $Z_\infty$-algebra structures are \emph{canonically} isomorphic to the trivial $Z_\infty$-algebra, the half-shuffle algebra.
\end{thm*}
Applying the symmetrisation functor, we obtain as a corollary:
\begin{cor*}
    Given a $B_\infty$-algebra on a vector space $V$ admitting a compatible $Z_\infty$-algebra structure, there is an explicit canonical isomorphism to the shuffle algebra.
\end{cor*}

Our goal is then to determine which $B_\infty$-algebras arise from $Z_\infty$-algebras. While this question of ``Zinbielification'' is a hard one, with ties to Rota-Baxter algebras \cite{Lemay2025} and divided powers structures \cite{DOKAS2009}, we choose to focus in on a particularly simple class of approach to ``Zinbielification'', which we call \emph{quasi-canonical}, and establish quasi-canonical families of $Z_\infty$-algebras symmetrising to classical quasi-shuffle algebras.

\vspace{2em}
We proceed as follows. In Section 1, we outline the essential background, defining Zinbiel algebras, coalgebras, comodules, and $B_\infty$-algebras. In Section 2, we define $Z_\infty$-algebras abstracts and establish a concrete description of them in terms of a bilinear map. Using this concrete description, we establish the existence and uniqueness of an alphabet-fixing isomorphism between a $Z_\infty$-algebra and the half-shuffle algebra.

In Section 3, we begin our exploration of the essential image of the symmetrisation function by defining a quasi-canonical $Z_\infty$-structure in terms of a fixed $B_\infty$-structure, and give a few examples of such structure. In Section 4, we focus on the properties of quasi-canonical structures in the context of quasi-shuffle algebras, describing completely all quasi-canonical structures for the quasi-shuffle algebras on the alphabet $\NN^n$.

In Section 5, we establish a relationship between quasi-canonical $Z_\infty$-algebra structures on quasi-shuffle algebras and semigroup cohomology \cite{novikov2001}, providing a theoretical framework for examples previously explored. Finally, in Section 6, we make some remarks on the computation of quasi-canonical structures and compute an example of a possible quasi-canonical structure on the block shuffle algebra.

\subsubsection*{Abusive terminology}
Throughout essentially all that follows, it will be necessary to fix a basis of a vector space $V$. As such, we will not refer to this fixed basis outside of initial definitions, but it should always be considered to be part of the data of a vector space. Furthermore, maps between such vector spaces $V\to W$ are assumed to map the basis of $V$ to a subset of the basis of $W$, up to rescaling. We allow ourselves this assumption as, with the exception of Definition \ref{def:quasi-canon}, fixing this basis mostly serves to simplify notation in our arguments, and plays no meaningful role in their content.

In the context of the tensor coalgebra on $V$, we refer to the basis as an \textit{alphabet}, elements of the basis as \textit{letters} and tensor products of these elements as \textit{words}

In the context of the symmetric algebra on $V$, we identify this with a (non-unital) polynomial algebra and use all associated terminology.

\section{Preliminaries}
\subsection{Zinbiel algebras}
\begin{definition}
    A Zinbiel algebra $Z$ over a field $K$ is a vector space $Z$ equipped with a bilinear product 
    \begin{align*}
        Z\otimes Z &\to Z\tx{,}\\
        a\otimes b &\mapsto a\succ b\tx{,}
    \end{align*}
    such that
    \begin{equation}\label{eq:Zinbiel-condition}
        (a\succ b)\succ c = a\succ (b\succ c + c\succ b)
    \end{equation}
    for all $a,b,c\in Z$.
\end{definition}

The structure of a Zinbiel algebra on a vector space $Z$ also endows $Z$ with the structure of a (non-unital) commutative associative algebra via the symmetrized product
\[a\star b := a\succ b + b\succ a.\]

    We define our Zinbiel algebras to be non-unital, as it is more convenient for our purposes to not have to make exceptions for the unit element in our statements. However, this is not always the case, and a Zinbiel product can be extended to one on $Z_1=K\oplus Z$ by defining
    \[ a\succ 1 = a,\quad 1\succ a = 0\tx{,}\]
    and restricting the condition given in Equation (\ref{eq:Zinbiel-condition}) to $a,b,c\in Z$. With this extension, we define the notion of a Zinbiel structure on an augmented algebra.
    
    Let $A$ be a commutative associative algebra over a field $K$ of characteristic $0$ with unit $\eta:K\to A$, product
    \begin{align*}
        A\otimes A &\to A\tx{,}\\
        a \otimes b &\mapsto a\star b\tx{,}
    \end{align*}
    and augmentation map $\varepsilon:A\to K$.

\begin{definition}
    A Zinbiel structure on an augmented algebra $A$ is a right $A$-module structure on $\overline{A}:=\ker \varepsilon$     
    \[\succ:\overline{A}\otimes A\to \overline{A}\]
    such that
    \[a\succ b + b\succ a =a\star b, \quad \text{and}\quad a\succ 1=1\]
    for all $a,b\in \overline{A}\subset A$.
\end{definition}

\begin{remark}
    Note that a Zinbiel structure on $A$ automatically endows $\overline{A}$ with the structure of a Zinbiel algebra.
\end{remark}

\subsection{Coalgebras, comodules, and \texorpdfstring{$B_\infty$}{B-infinity}-algebras}
At the core of this topic is the notion of a coalgebra and a comodule. In all that follows, fix a field $K$ of characteristic $0$.
\begin{definition}
    A coalgebra $(C,\Delta,\varepsilon)$ is a vector space over a field $K$ equipped with a coproduct
    \[\Delta:C\to C\otimes C\]
    and counit $\varepsilon: C\to K$ such that
    \begin{enumerateroman}
        \item $(\id_C\otimes \Delta)\circ \Delta = (\Delta\otimes\id_C)\circ\Delta$,
        \item $(\id_C\otimes \varepsilon)\circ \Delta = \id_C = (\varepsilon\otimes \id_C)\circ \Delta$.
    \end{enumerateroman}

    A coaugmentation of $C$ is a non-zero element $1\in \ker\varepsilon$ such that $\Delta 1 = 1\otimes 1$. For a coaugmented coalgebra, we defined the reduced coproduct
    \[\Delta_2(x) := \Delta(x) - x\otimes 1 - 1\otimes x\]
    and the iterated reduced coproduct
    \[\Delta_{n+1} = (\Delta\otimes \id_C^\otimes{n-1})\circ \Delta_n.\]
    Finally, we call $C$ conilpotent if, for every $x\in C$, there exists $n\geq 2$ such that $\Delta_n(x)=0$.
\end{definition}

\begin{definition}
    Given a coalgebra $C$ over a field $K$, a right-comodule $M$ over $C$ is a $K$-vector space with a linear map
    \[\rho:M\to M\otimes C\]
    such that
    \begin{enumerateroman}
        \item $(\id_M\otimes \Delta)\circ \rho=(\rho\otimes \id_C)\circ \rho$,
        \item $(\id_M\otimes \varepsilon)\circ \rho = \id_M$.
    \end{enumerateroman}
    If $C$ is coaugmented, we can define the reduced coaction by
    \[\rho_2(m) = \rho(m) - m\otimes 1\]
    and similarly define an iterated reduced coaction $\rho_n$. We call $M$ conilpotent if for all $m\in M$ there exists $n\geq 2$ such that $\rho_n(m)=0$.
\end{definition}

\begin{example}
    Let $V$ be a vector space over $K$. The tensor coalgebra $T^c(V)$ is the conilpotent coalgebra with underlying vector space
    \[T^c(V) = \bigoplus_{n\geq 0} V^{\otimes n},\]
    coaugmentation given by $1\in K=V^{\otimes 0}$, and coproduct given by deconcatenation
    \begin{align*} 
    \Delta(v_1\otimes v_2\otimes \cdots v_n) =&{} (v_1\otimes\cdots\otimes v_n)\bigotimes 1 + 1\bigotimes (v_1\otimes \cdots \otimes v_n)\\
    &+\sum_{i=1}^{n-1}(v_1\otimes \cdots \otimes v_i)\bigotimes (v_{i+1}\otimes \cdots \otimes v_n).
    \end{align*}
    This is the cofree conilpotent coalgebra on $V$ \cite{Sweedler1969-kb}, in that $T^c(V)$ is universal among all conilpotent coalgebras $C$ with a linear map $C\to V$.
\end{example}

In order to simplify notation, we will suppress the tensor product and write
\[v_1v_2\ldots v_n := v_1\otimes v_2\otimes\cdots\otimes v_n.\]
As mentioned earlier, we will refer to a tensor product of basis elements as a word in $V$

A $B_\infty$-algebra is an algebra structure on $T^c(V)$ for some $V$ that upgrades $T^c(V)$ to a Hopf algebra. Note that, as $T^c(V)$ is conilpotent, the notions of Hopf algebra and bialgebra are equivalent, and so we do not need to make reference to an antipode.

\begin{definition}
    A $B_\infty$-algebra on a vector space $V$ is a coalgebra morphism
    \[m:T^c(V)\otimes T^c(V) \to T^c(V)\]
    such that $(T^c(V),m,1)$ is an algebra and 
    \[\Delta:T^c(V)\to T^c(V)\otimes T^c(V)\]
    is an algebra morphism.
\end{definition}

As in \cite{binfty_foissy_patras}, a $B_\infty$-algebra on $V$ is determined uniquely by a bilinear map
\[(-,-):T^c(V)\otimes T^c(V)\to V\]
such that
\begin{enumerateroman}
    \item \[ (1,v_1v_2\ldots v_n) = (v_1v_2\ldots v_n,1) = \begin{cases} v_1 \text{ if }n=1,\\ 0\text{ otherwise.}\end{cases}\]
    \item For all non-empty words $u,v,w$ in $V$
    \[(u,v\star w) = (u\star v,w)\]
    where
    \[ u\star v := \sum_{k\geq 1}\sum_{\substack{u=u_1u_2\ldots u_k\\ v=v_1v_2\ldots v_k}} (u_1,v_1)(u_2,v_2)\ldots (u_k,v_k)\]
    is given as a sum over splittings of $u,v$ into subwords.
\end{enumerateroman}
We call the $B_\infty$-algebra commutative if $(u,v)=(v,u)$. 

\begin{theorem}[\cite{binfty_foissy_patras}]\label{theorem:foissy}
    Every commutative $B_\infty$-algebra on $V$ is isomorphic to the $B_\infty$-algebra with $(u,v)=0$ for all pairs of words such that $u\neq 1$ or $v\neq 1$.
\end{theorem}

\begin{remark}
Note also that the map associating a $B_\infty$-algebra to a vector space $V$ with an appropriate bracket is functorial in the following sense. Given $B_\infty$-algebras $C$ over $V$ and $D$ over $W$, and a linear map $\phi:V\to W$ such that the induced map $\phi:T^c(V)\to T^c(W)$ satisfies
\[\phi\left((a,b)_C\right) = \left(\phi(a),\phi(b)\right)_D\]
is is easy to verify that $\phi:C\to D$ is a map of $B_\infty$-algebras:
\begin{align*}
    \phi(u\star v) &= \sum_{n=1}^\infty \sum_{\substack{u=u_1\ldots u_n,\\ v=v_1\ldots v_n}}\phi\left((u_1,v_1)\ldots (u_n,v_n)\right)\\
    &= \sum_{n=1}^\infty \sum_{\substack{u=u_1\ldots u_n,\\ v=v_1\ldots v_n}} \left(\phi(u_1),\phi(v_1)\right)\ldots\left(\phi(u_n),\phi(v_n)\right)\\
    &= \sum_{\substack{\phi(u)=u^\prime_1\ldots u^\prime_n,\\ \phi(v)=v^\prime_1\ldots v^\prime_n}} (u^\prime_1,v^\prime_1)\ldots (u^\prime_n,v^\prime_n)\\
    &= \phi(u)\star \phi(v).
\end{align*}
\end{remark}

\section{\texorpdfstring{$Z_\infty$}{Z-infinity}-algebras and their classification}
Over a given vector space $V$, consider the cofree conilpotent coalgebra $C=T^c(V)$ and the category of conilpotent right comodules $M$ with a surjective map $M\to V$. We claim $\overline{C}:=T^c_+(V)$ equipped with coaction
\[\overline{\Delta}(x):=\Delta(x)-1\otimes x\]
and natural projection $\pi|_{\overline{C}}:\overline{C}\to V$ is an terminal object in this category, and hence merits the title of a cofree conilpotent comodule.

\begin{remark}
    As we will always be taking our structures to be conilpotent, we will no longer explicitly assume conilpotence of any coalgebra or comodule.
\end{remark}

\begin{proposition}
    Given a right $C$-comodule $M$ with a surjective map $\phi:M\to V$, there exists a unique map of comodules $\Phi:M\to \overline{C}$ such that $\pi|_{\overline{C}}\circ\Phi=\phi$.
\end{proposition}
\begin{proof}
    Consider the map 
    \[\Phi = \sum_{n=1}^\infty \phi\otimes\pi^{\otimes n-1}\Delta^n_M\]
    obtained by composing the iterated coaction with the given projection maps. By conilpotence of $C$, $\Phi(m)$ is a finite sum for every $m\in M$, and hence this is a well defined map of vector spaces
    \[\Phi:M\to \overline{C}.\]
    It is not hard to see that this is a map of comodules satisfying the given condition, so it remains to check that this is unique. Indeed, suppose we have such a $\Phi:M\to\overline{C}$. Then we can write
    \[\Phi= \sum_{n=1}^\infty \Phi_n\]
    using the natural grading on $\overline{C}$. By the definition of the coaction in $\overline{C}$,
    \begin{align*} 
    \Phi_n(m) &= \pi^{\otimes n} \overline{\Delta}^n (\Phi(m))\\
                &= \pi^{\otimes n} (\Phi\otimes\id^{\otimes n-1}\Delta^{n-1}_M(m)\\
                &= \phi\otimes\pi^{\otimes n-1}\Delta^{n-1}_M(m),
    \end{align*}
    from which the claim follows.
\end{proof}

Let us now assume that $C$ has been equipped with the structure of a fixed $B_\infty$-algebra, so that $C$ is a bialgebra. We then define a $Z_\infty$-algebra as a special Zinbiel structure on $C$.

\begin{definition}
    A $Z_\infty$-algebra on $(C,\star)$ is a Zinbiel structure on $C$
    \[\succ:\overline{C}\otimes C \to \overline{C}\]
    such that $\succ$ is a morphism of $C$-comodules, and $\overline{\Delta}$ is a morphism of $C$-modules.
\end{definition}

\begin{lemma}
    The structure of a $Z_\infty$-algebra on a $B_\infty$-algebra $(C,\star)$ is entirely encoded by the bilinear map
    \[\langle -,-\rangle:\overline{C}\otimes C \to V\]
    given by
    \[\langle u,v\rangle = \pi(u\succ v).\]
    The Zinbiel product is then given by
    \[u\succ v = \sum_{\substack{u=u_1u_2,\\v=v_1v_2}} \langle u_1,v_1\rangle (u_2\star v_2).\]
\end{lemma}
\begin{proof}
    We first note that the bracket $\langle -,-\rangle$ defines a surjective map $\overline{C}\otimes C\to V$, and hence there is a unique map of comodules 
    \[\Phi:\overline{C}\otimes C \to \overline{C}\]
    such that 
    \[\langle u,v\rangle = \pi(\Phi(u\otimes v))\]
    Furthermore, this map $\Phi$ is given by
    \begin{align*}
        \Phi(u,v) &=\sum_{n\geq 1} \sum_{\substack{u=u_1u_2\ldots u_n ,\\v=v_1v_2\ldots v_n}}\langle u_1,v_1\rangle\pi(u_2\star v_2)\cdots \pi(u_n\star v_n)\\
        &=\sum_{n\geq 1} \sum_{\substack{u=u_1u_2\ldots u_n\ldots ,\\v=v_1v_2\ldots v_n}}\langle u_1,v_1\rangle\pi(u_2\star v_2)\cdots \pi(u_n\star v_n)\\
        &=\sum_{\substack{u=u_1u_2,\\v=v_1v_2}} \langle u_1,v_1\rangle (u_2\star v_2)
    \end{align*}
\end{proof}

    Note that the $Z_\infty$-bracket captures not just the $Z_\infty$-algebra structure, but also the $B_\infty$-algebra structure, which is determined by the map
    \begin{align*}
        C\otimes C &\to C\\
        u\otimes v &\mapsto (u,v):= \langle u,v\rangle + \langle v,u\rangle
    \end{align*}
    where we define $\langle 1,v\rangle =0$ for all non-trivial words $v$.

\begin{example}\label{ex:trivial}
On any alphabet $X$, we can define a trivial $B_\infty$-algebra structure in which 
\[(u,v) = 0\]
unless $u=1$ or $v=1$. This leads to the classical shuffle algebra $(\bbQ\langle X\rangle, \shuffle)$. There is an associated trivial $Z_\infty$-algebra structure on $\bbQ\langle X\rangle$, determined by
\[\langle u,v\rangle = 0 \text{ for all }v\neq 1\]
which defines the free Zinbiel algebra on $X$, with product
\[au\succ bv = a(u\shuffle bv).\]
\end{example}

\begin{example}\label{ex:quasi-shuffle}
    Consider the alphabet $Y=\{Y_i\}_{i\in\bbZ_{>0}}$. We define a $B_\infty$-algebra structure on $\bbQ\langle Y\rangle$ by defining
    \[(y_i,y_j)=y_{i+j}\]
    and setting all other non-trivial brackets to zero. This leads to the harmonic algebra, a special case of a quasi-shuffle algebra. This admits a $Z_\infty$-algebra structure, given by
    \[\langle y_i,y_j\rangle = \frac{i}{i+j}y_{i+j}.\]
    So, for example:
    \begin{align*}
        y_i \star y_j &= y_iy_j + y_jy_i + y_{i+j},\\
        y_i \succ y_j &= y_iy_j + \frac{i}{i+j}y_{i+j}.
    \end{align*}
\end{example}

Given any bracket, we can define a product on $C$ recursively by $1\star 1:=1$ and
\begin{equation}\label{eq:binf-product}
    u\star v:= \sum_{\substack u=u_1u_2,\\v=v_1v_2} (u_1,v_1)(u_2\star v_2)
\end{equation}
and a $C$-module structure on $\overline{C}$ by
\begin{equation}\label{eq:zinf-product}
u\succ v :=\sum_{\substack u=u_1u_2,\\v=v_1v_2} \langle u_1,v_1\rangle (u_2\star v_2).
\end{equation}
Note that $u\star v = u\succ v + v\succ u$ for all $u,v$ with the above conventions. However, these do not necessarily define a $Z_\infty$-structure. 

Then we have the following
\begin{proposition}
    A bilinear bracket
    \[\langle -,-\rangle:\overline{C}\otimes C \to V\]
    defines a $Z_\infty$-algebra structure on $C$ if and only if
    \[\langle u\succ v,w\rangle = \langle u, v\star w\rangle\]
    for all $u,v,w\in \overline{C}$, using the products defined by Equations (\ref{eq:binf-product}) and (\ref{eq:zinf-product}).
\end{proposition}
\begin{proof}
        A short argument is given by noting that both
        \[\succ\circ(\id\otimes \star): \quad \text{and}\quad \succ\circ(\succ \otimes \id)\]
        define comodule morphisms from $\overline{C}^\otimes{2}\otimes C$ to $\overline{C}$ and are therefore entirely determined by their projection to $V$. We can also show this much more constructively.
        
    Clearly, if we have a $Z_\infty$-algebra structure on $C$, then
    \[(u\succ v)\succ w) = u\succ(v\star w)\]
    and hence the projection onto $V$ of each side must be equal, giving the desired identity among brackets. Now suppose that
    \[\langle u\succ v,w\rangle = \langle u, v\star w\rangle\]
    for all $u,v,w\in \overline{C}$. First note that this implies
    \begin{align*}
        (u\star v,w) &= \langle u\succ v, w\rangle> + \langle v\succ u,w\rangle + \langle w,u\star v\rangle\\
        &= \langle u,v\star w \rangle + \langle v,u\star w\rangle + \langle w\succ v, u\rangle\\
        &= \langle u, v\star w\rangle + \langle v\succ w,u\rangle + \langle w\succ v,u\rangle\\
        &= \langle u,v\star w\rangle + \langle v\star w,u\rangle\\
        &= (u,v\star w),
    \end{align*}
    using the commutativity of the $\star$-product.
    We then have
    \begin{align*}
        (u\succ v)&\succ w = \sum_{n\geq 1}\left(\sum_{\substack{u=u_1u_2\ldots u_n\ldots ,\\v=v_1v_2\ldots v_n}}\langle u_1,v_1\rangle(u_2,v_2)\cdots (u_n,v)n)\right)\succ w\\
        =&{} \sum_{n\geq 1} \sum_{i=1}^n \sum_{w=w_1w_2}\left(\sum_{\substack{u=u_1u_2\ldots u_n\ldots ,\\v=v_1v_2\ldots v_n}}\langle\langle u_1,v_1\rangle(u_2,v_2)\cdots (u_i,v_i),w_1\rangle (u_{i+1}v_{i+1})\ldots (u_n,v_n)\star w_2\right)\\
        =&{}\sum_{\substack{u=u_1u_2,\\v=v_1,v_2}} \sum_{w=w_1w_2} \langle u_1\succ v_1,w_1\rangle (u_2\star v_2)\star w_2\\
        =&{} \sum_{\substack{u=u_1u_2,\\v=v_1,v_2}}\sum_{w=w_1w_2}\langle u_1,v_1\star w_1\rangle u_2\star(v_2\star w_2)\\
        =&{} \sum_{n,i\geq 0}\sum_{u=u_1u_2} \sum_{\substack{v=v_1\ldots v_n,\\ w=w_1,\ldots,w_n}} \langle u_1, (v_1,w_1)\cdots (v_i,w_i)\rangle (v_{i+1},w_{i+1})\cdots (v_n,w_n)\\
        &= u\succ(v\star w)
    \end{align*}
\end{proof}

\subsection{Isomorphisms of \texorpdfstring{$Z_\infty$}{Z-infinity}-algebras}
As discussed previously, some algebras that are not finitely generated can be equipped with a Zinbiel structure such that they become finitely generated as Zinbiel algebras. As such, the set of morphisms compatible with a pair of given Zinbiel structures can be much simpler than the full set of homomorphisms. In the case of $Z_\infty$-algebras, they become extremely simple.

\begin{proposition}\label{prop:zinf-morphism}
    Given a $Z_\infty$-algebra structure $(C,\succ)$, the set of morphisms of Zinbiel algebras
    \[\phi:(C,\succ_{triv})\to (C,\succ)\]
    is in bijection with linear maps $V\to \overline{C}$. In particular, there is a unique Zinbiel morphism whose restriction to $V$ is the identity.
\end{proposition}
\begin{proof}
    Fix a basis $X=\{x_i\}_{i\in I}$ of $V$. The Zinbiel morphism induces a map $V\to \overline{C}$. Conversely, the set of words in $X$ forms a basis for $\overline{C}$, and we have that
    \[w_1w_2\ldots w_n = w_1\succ_{triv}(w_2\succ_{triv}(\cdots w_{n-1}\succ_{triv} w_n)\cdots )\]
    for every word $w_1\ldots w_n$. Thus, the image of $w_1\ldots w_n$ for any word is entirely determined by the images of $w_1,\ldots,w_n\in X$.
\end{proof}

\begin{remark}
    We sometimes refer to morphisms between algebras
    \[\phi:(C,\star_1)\to (C,\star_2)\]
    such that $\phi|_V=\id$ as begin \textit{tangent to the identity}. Proposition \ref{prop:zinf-morphism} implies that, given a pair $Z_\infty$-algebra structures, there is a unique map $\phi:C\to C$ that is a Zinbiel homomorphism and is tangent to the identity.
\end{remark}

We claim that the map induced by the identity on $V$ (or any automorphism of $V$), induces an isomorphism of Zinbiel algebras.

\begin{theorem}
    All $Z_\infty$-algebras $(C,\succ)$ are isomorphic as Zinbiel algebras to the free Zinbiel algebra on $V$.
\end{theorem}
\begin{proof}
    As noted before, the identity induces a map of Zinbiel algebras, which we claim is surjective. To see this, note that we can filter $C$ by the length of words in a basis $X$ of $V$, and that the induced map is a map of filtered Zinbiel algebras. Furthermore
    \[\phi(w) = w + \text{ words of lower length.}\]
    Thus, the length graded pieces of this map are bijective, and hence the map itself is bijective, and hence has a inverse map, at least as a map of vector spaces. To see that the inverse is a Zinbiel map, note that
    \begin{align*}
        &{}\phi(\phi^{-1}(u\succ v) - \phi^{-1}(u)\succ_{triv}\phi^{-1}(v))\\
        =&{} u\succ v - \phi(\phi^{-1}(u))\succ\phi(\phi^{-1}(v))\\
        =&0.
    \end{align*}
\end{proof}

As an immediate corollary, we obtain the following special case of Theorem \ref{theorem:foissy} \cite{binfty_foissy_patras}, giving a much simpler argument in this special case.

\begin{corollary}
    Every $B_\infty$-algebra that admits a $Z_\infty$-structure is isomorphic to the shuffle algebra.
\end{corollary}

Thus, the question is then which commutative $B_\infty$-algebras admit a $Z_\infty$-structure? 

\section{\texorpdfstring{$Z_\infty$}{Z-infinity}-algebras associated to commutative \texorpdfstring{$B_\infty$}{B-infinity}-algebras}
 It is necessary that a $B_\infty$-algebra they be commutative in order to admit a $Z_\infty$-structure. However, we can find no clear condition beyond this. Nevertheless, we make the following (optimistic) conjecture

\begin{conjecture}
    Every commutative $B_\infty$-algebra containing no basis element such that
    \[(v,v^\prime)=v\text{ for all }v\in V\]
    admits a $Z_\infty$-structure
\end{conjecture}

A naive method to construct a $Z_\infty$-structure from a $B_\infty$-structure is to demand that the associated brackets be proportional.

\begin{definition}\label{def:quasi-canon}
    Given a $B_\infty$-algebra over an alphabet $X$, we call an admissible $Z_\infty$-structure quasi-canonical if for every pair of words $u,v$, we have
    \[\langle u,v\rangle = \eta_{u,v}(u,v)\]
    for some function $\eta:X^\times\times X^\times \to K^\times$.
\end{definition}

\begin{remark}
    This definition is the first to explicitly rely on a choice of basis $X$ of the underlying vector space $V$.
\end{remark}

\begin{example}
    Consider the quasi-shuffle algebra \cite{hoffman00} $C=\bbQ\langle Y\rangle$ on the alphabet 
    \[Y=\{y_1,y_2,y_3,\ldots\}.\]
    This is a commutative $B_\infty$-algebra, with bracket given by
    \[(y_i,y_j)=y_{i,j},\quad (u,1)=(1,u)=u\]
    and $(u,v)=0$ for all other pairs of words $u,v$. As seen in Example \ref{ex:quasi-shuffle}, this admits a quasi-canonical $Z_\infty$-structure, determined by
    \[\langle y_i,y_j\rangle = \frac{i}{i+j}y_{i+j}\]
    As will be shown in Proposition \ref{prop:quasi-shuffle}, this quasi-canonical structure is in fact unique. This is not generally the case, as illustrated by the next example.
\end{example}

\begin{example}
   Consider the $B_\infty$-algebra $\bbQ\langle Y\rangle$ on the alphabet
   \[Y=\{y_{i,j}\}_{\substack{i,j\in \ZZ_{\geq 0}\\ (i,j)\neq (0,0)}}\]
   determined by the bracket
   \[(y_{i,j},y_{k,\ell}) = y_{i+k,j+\ell},\quad (u,1)=(1,u)=u,\]
   and $(u,v)=0$ for all other pairs of words $(u,v)$. Then, for any $\alpha\in \QQ$, the bracket
   \[\langle y_{i,j},y_{k,\ell}\rangle = \frac{i+\alpha j}{i+k +\alpha(j+\ell)}y_{i+j,k+\ell}\]
   defines a quasi-canonical $Z_\infty$-structure.
\end{example}
Note that a quasi-canonical structure does not necessarily exist.

\begin{example}\label{ex:no-quasi-canon}
    Consider the $B_\infty$-algebra over the vector space $V=\bbQ a\oplus \bbQ b$, determined by the bracket
    \[(a,a)=b,\quad (a,b)=b,\quad (b,b)=b.\]
    If this admits a quasi-canonical $Z_\infty$-structure, then we must have that
    \[a\succ a= a^2 + \frac{1}{2}b,\quad \text{ and }\quad a\succ b = ab + \alpha b\]
    for some $\alpha\in \bbQ$. Comparing
    \[(a\succ a)\succ a = a\succ (a\star a)\]
    we find $\alpha=\frac{1}{3}$. Next consider
    \[\left((a\succ a)\succ a\right)\succ a = a\succ( a\star a\star a).\]
    Computing the lowest length term on each side, we find that we must have
    \[\frac{2}{9}b = \frac{1}{3}b\]
    which is impossible. A similar issue arises in any $B_\infty$-algebra for which there exists a ``zero'' element $b$ such that $(z,b)=b$ for all letters $z$.
\end{example}

The question of whether a quasi-canonical structure exists is in general hard, but in many cases it is possible to describe a finite family of potential structures, as illustrated in Theorem \ref{theorem:blk-shuffle}. If the associated $B_\infty$-algebra is simple enough, we can much more explicitly explore quasi-canonical structures. 

\begin{remark}
    Introducing quasi-canonical structures serves to provide a comparatively simple approach to understanding if we can Zinbielify an commutative algebra. However, this is not the only perspective one could take. For example, one could approach this questions via the equivalence of Zinbiel algebras and FTC-algebras \cite{Lemay2025}, giving ties to the theory of Rota Baxter algebras; or using the existence of a canonical divided powers structure in characteristic $0$ to construct an adjoint functor to symmetrisation \cite{DOKAS2009}.
\end{remark}

\section{Quasi-canonical structures for quasi-shuffle algebras}
A quasi-shuffle algebra is usually defined recursively as follows. Given an alphabet $A$, and a commutative, associative map
\[(-,-):A\times A\to KA\]
we define a product on $K\langle A\rangle$ recursively by
\begin{align*}
u\star 1 = 1 \star u = u \text{ for all u,}\\
au\star bv = a(u\star bv) + b(au\star v) + (a,b)(u\star v)
\end{align*}
for all letters $a,b\in A$ and words $u,v$. We can view this as a special case of a $B_\infty$-algebra as follows.

\begin{definition}
    A quasi-shuffle algebra $(C,\star)$ over a vector space $V$ is a $B_\infty$-algebra with bracket such that $(u,v)=0$ unless $u=1$, $v=1$ or $|u|=|v|=1$.
\end{definition}

The bracket induces an algebra structure on $V$, and so we can view a quasi-shuffle algebra as a structure associated to an arbitrary non-unital algebra. When studying the existence of quasi-canonical $Z_\infty$-structures, these provide an excellent source of examples, as $Z_\infty$-conditions are particularly simple.

\begin{lemma}\label{lem:quasi-shuffle-bracket}
    A bracket on a quasi-shuffle algebra $C$ over $V$ satisfying
    \[\langle \langle a,b\rangle,c\rangle = \langle a, (b,c)\rangle\]
    for all non-trivial $a,b,c\in V$ induces a $Z_\infty$-structure on $C$.
\end{lemma}
\begin{proof}
    It is clear that we can define $\langle u,v\rangle=0$ if $(u,v)=0$ without issue. As such, the $Z_\infty$-condition
    \[\langle u\succ v,w\rangle = \langle u, v\star w\rangle\]
    is trivial unless $|u|=|w|=1$. The only term in $u\succ v$ that can contribute is the term of length $1$, given by $\langle u,v\rangle$, which is $0$ unless $|v|=1$. Similarly, the only term that can contribute in $v\star w$ is the term of length $1$ $(v,w)$, which is $0$, unless $|v|=1$. Thus, the condition is trivial unless $u,v,w$ are all words of length $1$, in which case the condition reduces to
    \[\langle \langle u,v\rangle,w\rangle = \langle u, (v,w)\rangle.\]
\end{proof}

Every quasi-canonical structure on a quasi-shuffle algebra must come from such a bracket. We can easily construct quasi-canonical structures for quasi-shuffle algebras associated to $\Sym(V)$, viewed as a non-unital algebra, for finite dimensional $V$.

\begin{proposition}\label{prop:quasi-shuffle}
The quasi-shuffle algebra $C=\QQ\langle Y\rangle$ on the alphabet $Y=\{y_1,y_2,\ldots\}$ determined by
\[(y_i,y_j)=y_{i+j},\quad (y_i,1)=(1,y_i)=y_i\]
admits a unique quasi-canonical $Z_\infty$-structure.
\end{proposition}
    \begin{proof}
A quasi-canonical $Z_\infty$-algebra structure is entirely determined by the bracket
\[\langle y_i,y_j\rangle = \eta_{i,j}y_{i+j}\]
In particular, we must have that
    \[y_i\succ y_j = y_iy_j + \eta_{i,j}y_{i+j}.\]
    The Zinbiel condition then implies that
    \[\eta_{i,j}\langle y_{i+j},y_k\rangle = \langle y_i, y_j\star y_k\rangle\]
    and hence
    \[\eta_{i,j}\eta_{i+j,k} = \eta_{i,j+k}\]
    for all $i,j,k$. We must also have that $\eta_{i,j}+\eta_{j,i}=1$, and so $\eta_{i,i}=\frac{1}{2}$. We next have that
    \[\eta_{1,j}\eta_{1+j,1} =\eta_{1,1+j} = 1-\eta_{1+j,1}\]
    which we can solve inductively to find that
    \[\eta_{1,n}=\frac{1}{n+1},\quad \eta_{n,1}=\frac{n}{n+1}\]
    Finally, we have that
    \[\eta_{1,i-1}\eta_{i,j}=\eta_{1,i+j-1}\]
    and so
    \[\eta_{i,j}=\frac{i}{i+j}.\]
    It is then easy to check that this does indeed satisfy the Zinbiel condition, and so we obtain a $Z_\infty$-algebra structure. 
\end{proof}

    We therefore obtain a quasi-canonical isomorphism of Zinbiel algebras, and hence of algebras, between the shuffle algebra and the quasi-shuffle algebra, that is distinct from those of Hoffman \cite{hoffman00}, and Foissey and Patras \cite{binfty_foissy_patras}. In low degree, it is given by
    \[\phi(y_n)=y_n,\quad \phi(y_my_n)=y_my_n + \frac{m}{m+n}y_{m+n}\]
    and so on.

This can be generalised easily to a quasi-canonical $Z_\infty$-structure on a multi-index analogue of the given quasi-shuffle algebra. This generalisation of the quasi-shuffle algebra is important enough that we will give it a name

\begin{definition}
    Denote by $\bbN^n$ the free abelian semigroup generated by $\{e_1,e_2,\ldots,e_n\}$, and let $Y_n=\{y_{\vec{m}}\mid \vec{m}\in\bbN^n\}$ be the alphabet indexed by $\bbN^n$. Define the $n$-quasi-shuffle algebra $C_n$ to be the quasi-shuffle algebra on $\QQ\langle Y\rangle$ to be determined by
    \[(y_{\vec{m}_1},y_{\vec{m}_2}) = y_{\vec{m}_1+\vec{m}_2},\quad (y_{\vec{m}},1)=(1,y_{\vec{m}})=y_{\vec{m}}\]
\end{definition}
\begin{proposition}\label{prop:multi-index}
    The $n$-quasi-shuffle algebra $C_n$ admits a quasi-canonical $Z_\infty$-structure for every tuple $(c_1,\ldots,c_{n-1})\in\QQ_{>0}^{n-1}$.
\end{proposition}
\begin{proof}
    Suppose we know $\eta_{e_i,e_j}$ for each $i\neq j$. Let $\vec{m}\in\bbN^n\setminus\{e_1,\ldots,e_n\}$ and let $a$ be the smallest index for which $\vec{m}-e_a\in\bbN^n$. Then
    \[\eta_{e_a,\vec{m}-e_a}\eta_{\vec{m},\vec{m}^\prime} = \eta_{e_a,\vec{m}+\vec{m}^\prime-e_a}\]
    for any $\vec{m}^\prime\in \bbN^n$. Thus, every $\eta_{\vec{m},\vec{m}^\prime}$ is can be determined in terms of
    \[\eta_{e_a,\vec{m}}\text{ for }1\leq a \leq n,\; \vec{m}\in\bbN^n\]
    Using that
    \[\eta_{e_a,\vec{m}}\eta_{e_a+\vec{m},e_a} = \eta_{e_a,\vec{m}+e_a} = 1-\eta_{e_a+\vec{m},e_a}\]
    we can inductively determine that $\eta_{e_a,\vec{m}}$ in terms of the $\eta_{e_i,e_,j}$.
    
     Furthermore, we have that
    \begin{align*}
        \eta_{e_1,e_i}\eta_{e_1+e_i,e_j} = &\eta_{e_1,e_i+e_j} = \eta_{e_1,e_j}\eta_{e_1+e_j,e_i},\\
        \eta_{e_i,e_1}\eta_{e_1+e_i,e_j} = &\eta_{e_i,e_1+e_j} = \eta_{e_i,e_j}\eta_{e_i+e_j,e_1},\\
        \eta_{e_j,e_1}\eta_{e_1+e_j,e_i} = &\eta_{e_j,e_1+e_i} = \eta_{e_j,e_i}\eta_{e_i+e_j,e_1}.
    \end{align*}
    Rearranging these, we find
    \[\eta_{e_1,e_i}\eta_{e_i,e_j}\eta_{e_j,e_1} = \eta_{e_i,e_1}\eta_{e_j,e_i}\eta_{e_1,e_j}\]
    Then, as $\eta_{\vec{m}_1,\vec{m}_2}+\eta_{\vec{m}_2,\vec{m}_1}=1$, we therefore have that 
    \[\eta_{e_i,e_j} = \frac{\eta_{e_i,e_1}\eta_{e_1,e_j}}{\eta_{e_i,e_1}\eta_{e_1,e_j}+\eta_{e_j,e_1}\eta_{e_1,e_i}}\]
    is uniquely determined by $\eta_{e_1,e_i}$ and $\eta_{e_1,e_j}$. It is straightforward to verify that this is consistent with all relations among $\{\eta_{e_i,e_j}\}$. Setting $\eta_{1,e_i}=c_{i-1}$, the claim follows.
\end{proof}

\begin{remark}
    We require that the tuple $(c_1,\ldots,c_{n-1})$ consist of positive rationals to avoid division by $0$. If, for example, $n=2$ and $c_1=-\frac{r}{s}$ for some $r,s\in\bbN$, then one can check that $\eta_{(r+s)e_1,re_2}$ is not well defined. The results of Section \ref{sec:qc-from-coho} help clarify this.
\end{remark}

\subsection{Examples from functoriality}
As noted in the preliminaries, the map taking a pair $(V,(-,-))$ to the associated $B_\infty$-algebra is functorial in an appropriate sense. In the case of quasi-shuffle algebras, we obtain a genuine functor
\[\QS:\NAlg_K \to \GrAlg_K\]
from non-unital algebras over $K$ to graded algebras over $K$.

Given a non-unital algebra $A$ such that $\QS(A)$ admits a $Z_\infty$-structure for which $\langle u,v\rangle =0$ unless $|u|=|v|=1$, the bracket
\[\langle -,-\rangle :A\otimes A\to A\]
defines a Zinbiel structure on $A$.

Conversely, given a Zinbiel structure on $A$, we obtain a $Z_\infty$-structure on $\QS(A)$ via Lemma \ref{lem:quasi-shuffle-bracket}. Thus, we have the following.

\begin{proposition}
    There is a bijection between $Z_\infty$-structures on $\QS(A)$ for which  $\langle u,v\rangle =0$ unless $|u|=|v|=1$, and Zinbiel structures on $A$.
\end{proposition}

As in the preliminaries, it is not difficult to see that the map taking a Zinbiel structure on $A$ to a $Z_\infty$-structure on $\QS(A)$ is functorial with respect to morphisms $A\to B$ of Zinbiel algebras. With this in mind, we can construct many examples of quasi-canonical $Z_\infty$-algebra structures.

\begin{proposition}\label{prop:ideal-mapping}
        Let $V$ be an $n$-dimensional $K$-vector space, and let $C=\QS\left( \Sym(V)\right)$ be the $n$-quasi-shuffle algebra $C_n$ on a basis of $V$. Let $A$ be a non-unital commutative algebra generated by $V$, so that there exists a unique algebra map
        \[\phi:\Sym(V)\to A\]
        whose restriction to $V$ can be taken to be the identity.

        Then, $\QS(A)$ admits a quasi-canonical $Z_\infty$-structure if ~$\ker(\phi)$ is a Zinbiel ideal for a quasi-canonical Zinbiel structure on $\Sym(V)$.
\end{proposition}
\begin{proof}
    Fix some choice quasi-canonical Zinbiel structure on $\Sym(V)$, by which we mean a Zinbiel structure such that
    \[a\succ_A b = \eta_{a,b}a\cdot b\]
    for some scalar $\eta_{a,b}$. Suppose $\ker(\phi)$ is a Zinbiel ideal. Then, we define, for any $a,b\in A$,
    \[a\succ_A b = \phi(\langle v_a,v_b\rangle)\]
    for any choice of lifts $v_a,v_b$ of $a$ and $b$ to $\Sym(V)$ respectively. This map is well defined: suppose $w_a,w_b$ are different lifts of $a$ and $b$. Then $v_a-w_a,v_b-w_b\in \ker(\phi)$. As $\ker(\phi)$ is a (two-sided) Zinbiel ideal, we therefore have
    \[\langle v_a,v_b\rangle -\langle w_a,w_b\rangle =\langle v_a-w_a,v_b\rangle + \langle w_a, v_b-w_b\rangle \in \ker(\phi).\]

    It is similarly straightforward to check that this defines a Zinbiel-structure on $A$, and furthermore that it is quasi-canonical. Hence, we obtain a quasi-canonical structure on $\QS(A)$.  
\end{proof}

\begin{example}
    Consider the quasi-shuffle algebra on the alphabet $Y=\{y_1,y_2\}$ with $B_\infty$-bracket given by
    \[(y_1,y_1)= y_2,\quad (y_1,y_2)=0,\quad (y_2,y_2)=0.\]
    This is the quasi-shuffle algebra associated to the non-unital algebra $x\bbQ[x]/(x^3)$, and it is clear to see that 
    \[(x^3)= \ker\left(x\bbQ[x]\to x\bbQ[x]/(x^3)\right)\]
    is a Zinbiel ideal for the unique quasi-canonical Zinbiel-structure on $x\bbQ[x]$ coming from Proposition \ref{prop:quasi-shuffle}. Thus, the bracket
    \[\langle y_1,y_1\rangle= \frac{1}{2}y_2,\quad \langle y_1,y_2\rangle=0,\quad \langle y_2,y_2\rangle =0\]
    defines a quasi-canonical $Z_\infty$-structure.
\end{example}

A partial converse holds, given conditions on the ideal $\ker(\phi)$.

\begin{proposition}\label{prop:ideal-converse}
    Let $V$ be an $n$-dimensional $K$-vector space, and let $C=\QS\left( \Sym(V)\right)$ be the $n$-quasi-shuffle algebra $C_n$ on a basis of $V$. Let $A$ be a non-unital commutative algebra generated by $V$, so that there exists a unique algebra map
        \[\phi:\Sym(V)\to A\]
        whose restriction to $V$ can be taken to be the identity.

        Then, if  $\QS(A)$ admits a quasi-canonical $Z_\infty$-structure and  $\ker(\phi)$ admits a generating set consisting of non-monomial elements of $\Sym(V)$, $\ker(\phi)$ is a Zinbiel ideal for a quasi-canonical Zinbiel structure on $\Sym(V)$.
\end{proposition}
\begin{proof}
       Use the fixed basis of $V$ to identify $\Sym(V)$ with the non-unital polynomial algebra in $n$-variables. Suppose that $\QS(A)$ admits a quasi-canonical $Z_\infty$-algebra structure. In order to show that $\ker(\phi)$ is a Zinbiel ideal for some Zinbiel-structure on $\Sym(V)$, we will use an inductive argument on the (finite) number of generators of $\ker(\phi)$.
    
    Let $A^\prime$ be another non-unital algebra generated by $V$, with a map $\psi:A^\prime\to A$ such that the composition
    \[\Sym(V)\longrightarrow A^\prime \longrightarrow A\]
    is equal to $\phi:\Sym(V)\to A$. To simplify notation, we will represent the image of monomials $a\in\Sym(V)$ in both $A$ and $A^\prime$ also by $a$.
    
    Suppose further that $\ker(\psi)$ is a principal ideal. If $\ker(\psi)$ is not generated by a monomial, then
    \[a\succ_A b = \eta_{a,b}a\cdot b \neq 0\]
    for any monomials $a,b\in \Sym(V)$. Hence, we can define a bilinear map on $A^\prime$ by
    \[a\succ_{A^\prime} b = \eta_{a,b} a\cdot b\]
    for every pair of monomial $a,b$. By definition, this must define a quasi-canonical Zinbiel-structure on $A^\prime$, as it satisfies the necessary condition for monomials, and the conditions are bilinear. With this structure, $\psi$ clearly defines a map of Zinbiel algebras, and so $\ker(\psi)$ is a Zinbiel ideal.

    By our assumptions on $\ker(\phi)$, there exist non-monomials $f_1,\ldots,f_k\in \Sym(V)$ such that $\ker(\phi)=(f_1,\ldots,f_k)$. By iterating the above argument along the sequence
    \[\Sym(V)\to \Sym(V)/(f_1) \to \Sym(V)/(f_1,f_2)\to\cdots \to \Sym(V)/(f_1,\ldots,f_k)\cong A\]
     we obtain a quasi-canonical Zinbiel-structure on $\Sym(V)$ for which $\ker(\phi)$ is a Zinbiel ideal.
\end{proof}

\begin{example}
    Consider again Example \ref{ex:no-quasi-canon}. Here we considered the quasi-shuffle algebra associated to $x\bbQ[x]/(x^3-x^2)$. In $x\bbQ[x]$, we have that
    \[x\succ (x^3-x^2) = \frac{1}{4}x^4 - \frac{1}{3}x^3 \not\in (x^3-x^2)\]
    so we do not obtain a Zinbiel ideal, and hence cannot obtain a quasi-canonical $Z_\infty$-structure.
\end{example}

We expect this converse statement to hold more generally. Indeed, given a quasi-canonical structure on $\QS(A)$, it is possible to determine equations for nearly every constant $\eta_{a,b}$ to define a quasi-canonical structure on $\Sym(V)$ in terms of that on $\QS(A)$, possibly with a few free parameters. While the ``obvious'' equations are all compatible, and define the same constants, it is not clear that all relations among the constants $\eta_{a,b}$ will be satisfied by this lift.

\section{Semigroup alphabets and a cohomological obstruction}
The condition that 
\[\eta_{\vec{a},\vec{b}}\eta_{\vec{a}+\vec{b},\vec{c}} = \eta_{\vec{a},\vec{b}+\vec{c}}\]
is very reminiscent of a (multiplicative) cocycle condition in group cohomology or Hochschild cohomology \cite{hochschild1945cohomology}. Indeed, we can re-interpret the results of the previous section entirely within the context of semi-group cohomology.

\subsection{Semigroup cohomology}
For semigroups, there are a number of related, but not equivalent, notions of cohomology, with two of the main approaches coming from the work of Novikov \cite{novikov2001}, and Leech and Grillet\cite{grillet95}, in the case of commutative semigroups. While our problem is most easily interpreted in the formalism of Novikov, we can relate the two in the case in which we are interested in order to exploit the results of Grillet.

\begin{definition}
    Let $S$ be a semigroup, and $M$ an $S$-module. Let $C^n(S,M)=\ZZ\Hom(S^n,M)$ be the set of $n$-chains of $S$ with coefficients in $M$ for all $n\geq 0$. We define a map
    \[\partial:C^n(S,M)\to C^{n+1}(S,M)\]
    by
    \begin{align*}
    \left(\partial f\right)(s_1,\ldots,s_{n+1}) =&{} s_1\cdot f(s_2,\ldots, s_{n+1})\\
    &+ \sum_{i=1}^n (-1)^i f(s_1,\ldots,s_{i-1},s_is_{i+1},s_{i+2},\ldots,s_{n+1}) +(-1)^{n+1}f(s_1,\ldots,s_n).
    \end{align*}
    One can check that $\partial^2=0$, equipping the collection $\{C^n(S,M)\}$ with the structure of a cochain complex. We define the Novikov cohomology of $S$ to be the cohomology of this complex, and denote the $n$\textsuperscript{th} cohomology group by $H^n_N(S,M)$.
\end{definition}

\begin{remark}
Novikov cohomology, as defined, is an explicit presentation of Eilenberg-Maclane cohomology \cite{kostin01}, and is isomorphic to the Ext-group
\[H^n_N(S,M)\cong \Ext^n_{\bbZ[S]-\Mod}(\bbZ,M)\]
where $\bbZ$ is the trivial $S$-module.

In contrast, the definition of Leech cohomology for commutative semigroups is substantially more involved, and so we will not provide a definition here. Leech cohomology $H^n_L(S,A)$ depends not on an $S$-module $A$, but on an abelian group valued functor $A$ from the Leech category of $S$. However, in the case where $A$ is the constant functor $A(s)=M$ for all $s\in S$, Leech cohomology can also be identified with extensions in the category of $S$-modules \cite{Grillet2022}
\[H^n_L(S,M)\cong \Ext^n_{\bbZ[S]-\Mod}(\bbZ,M).\]
As such, we can identify the two theories for our purposes.
\end{remark}

\subsection{Quasi-canonical \texorpdfstring{$Z_\infty$}{Z-infinty}-structures from cohomology}\label{sec:qc-from-coho}
Suppose that we have an alphabet $Y$ indexed by a commutative non-unital semigroup $S$. Then the product
\[(y_{s},y_{t}) = y_{st},\quad (y_s,1)=(1,y_s)=1\]
defines a quasi-shuffle algebra on $\QQ\langle Y\rangle$. A quasi-canonical $Z_\infty$-structure corresponds to a function
\[\eta: S\times S\to \QQ^\times\]
such that
\[\eta_{st,u}\eta^{-1}_{s,tu}\eta_{st}=1, \quad \eta_{s,t}+\eta_{t,s}=1\]
Equivalently, this is a cocycle valued in the $S$-module $\QQ^\times$, where $\QQ^\times$ is equipped with the trivial action
\[s\cdot q = 1\]
for all $s\in S$, in Novikov's theory of semigroup cohomology \cite{novikov2001}. Equivalently, we can consider (non-unital) Hochschild cohomology, with cocycles valued in the $\ZZ[S]$-bimodule $\QQ^\times$, with actions
\[s\cdot q = 1,\quad q\cdot s = q.\]
If the semigroup cohomology $H^2(S,\QQ^\times)=0$, then every cocycle is also a coboundary, and so there would exist $g: S\to\QQ^\times$ such that
\[\eta_{s,t} = g^{-1}_{st}g_s\].
The condition $\eta_{s,t}+\eta_{t,s}=1$ then implies that
\[g_{st}=g_s+g_t.\]
So, for semigroups with vanishing cohomology, quasi-canonical structures are determined by homomorphisms $g:S\to (\QQ,+)$, up to multiplication by an element of $\QQ^\times$.

\begin{example}
    If $H^2(\bbN^n,\QQ^\times)=0$, then a quasicanonical structure is determined by a tuple of $n$ non-zero rational numbers $(g(e_1),g(e_2),\ldots,g(e_n))$, up to multiplication by an element of $\QQ^\times$. As such we can assume $g(e_1)=1$, and we obtain exactly the results of Proposition \ref{prop:multi-index}.
\end{example}
From the work of Grillet \cite{grillet95}, we understand the cohomology of free semigroups very well.

\begin{theorem}[Grillet, \cite{grillet95}]
    Let $S=\bbN^n$ be the free commutative semi-group on $n$-generators, and let $\cH(S)$ be the Leech category of $S$. Then, for any abelian group valued functor $A$ on $\cH(S)$, $H^n_L(S,A)=0$ for all $n\geq 2$.
\end{theorem}

\begin{corollary}
    For any $\bbN^n$-module $M$, $\rmH^2_N(\bbN^n,M)=0$. 
\end{corollary}

Thus, we obtain another proof of Proposition \ref{prop:multi-index}: every homomorphism $g:\bbN^n\to (\bbQ_{>0},+)$, up to the multiplicative action of $\bbQ_{>0}$ on itself, determines a quasi-canonical structures on the multi-index quasi-shuffle algebra with 
\[\eta_{\vec{a},\vec{b}} = \frac{g(\vec{a})}{g(\vec{a})+g(\vec{b})}.\]
Conversely, every question-canonical structure determines a cocycle, and hence a homomorphism $g:\bbN^n\to (\bbQ,+)$, but to a multiplicative action of $\bbQ^\times$.

\begin{remark}
    This relationship severely limits the possible semigroups underlying algebras admitting quasi-canonical structures. For example, consider the cyclic semigroup 
    \[S=\langle a\mid a^3=a\rangle\]
    This has vanishing second cohomology, and so $\eta_{a,b}$ must arise from a non-zero homomorphism $f:S\to (\bbQ,+)$. However no such homomorphism exists, and as we saw in Example \ref{ex:no-quasi-canon}, no quasi-canonical structure can exist.
\end{remark}

\begin{remark}
    Unlike many classical appearances of cocycles in classification problems, we don't get a bijection between isomorphism classes of quasi-canonical structures and cohomology classes. Indeed, suppose $\eta_{a,b}^{(1)}$ and $\eta_{a,b}^{2}$ are cocycles defining Zinbiel structures on $\QQ\langle Y\rangle$. If $\eta_{a,b}^{(1)}$ and $\eta_{a,b}^{(2)}$ are cohomologous, then $(\QQ\langle Y\rangle,\succ_{(1)}$ and $(\QQ\langle Y\rangle),\succ{(2)}$ are, in general, only isomorphic as $\QQ\langle Y\rangle$ modules, \emph{not} as Zinbiel algebras.
\end{remark}

\section{Constructing potential quasi-canonical structures}
Unlike the case of quasi-shuffle algebras, there is no clear cohomological interpretation of a quasi-canonical $Z_\infty$-structure for a general $B_\infty$-algebra. Nevertheless, there exists a somewhat general framework for constructing potential quasi-canonical structures for a given explicit example.

\begin{lemma}\label{lem:rewrite-to-one}
    Given a $B_\infty$-algebra over $V$, any $Z_\infty$-structure is entirely determined by the brackets $\langle z,w\rangle$, for letters $z\in V$ and arbitrary words $w$.
\end{lemma}
\begin{proof}
    We claim that any bracket $\langle u,v\rangle$ can be expressed in terms of brackets of the form $\langle z,w\rangle $. Clearly this is true for words $u$ of length $|u|=1$, so assume it holds true for all words $u$ of length $|u|\leq n$ for some $n$. Then, for a word $zu$ of length $n+1$, note that
    \[zu = z\succ u + g(z,u)\]
    where $g(z,u)$ is a linear combination of words of length at most $n$. Hence
    \begin{align*}
        \langle zu, v\rangle &= \langle z\succ u,v\rangle - \langle g(z,u),v\rangle\\
        &= \langle z,u\star v\rangle -\langle g(z,u),v\rangle
    \end{align*}
    which can be expressed in terms of brackets of the desired form by induction. The claim then follows.
\end{proof}

The goal is then to use known properties of the $B_\infty$-algebra structure to determine $\langle z,w\rangle$ in terms of a finite set of such expressions. While this approach cannot ensure the existence of a quasi-canonical structure, it does often produce a finite dimensional space of potential such structures, as we will illustrate in the case of the generalised quasi-shuffle algebra of \cite{keilthy22}.

\begin{theorem}\label{theorem:blk-shuffle}
   The $B_\infty$-algebra 
   \[C=\bbQ\langle z_1,z_2,z_3,\ldots\rangle\]
   defined by the bracket
    \[(z_{i_1}y_{i_2}\ldots z_{i_k},z_{j_1}z_{j_2}\ldots z_{j_\ell})=\begin{cases} (-1)^{\min(k,l)}z_{i_1+\cdots+i_k+j_1+\cdots+j_\ell},\text{ if }|k-\ell|=1\\ 0 \text{ otherwise}\end{cases}\]
    admits an at most one quasi-canonical $Z_\infty$-structure.
\end{theorem}
\begin{proof}
    By Lemma \ref{lem:rewrite-to-one}, it suffices to determine the brackets $\langle z_a,w\rangle$ for all $a\geq 1$ and words $w$. In order for the bracket $\langle u,v\rangle $ to define a quasi-canonical structure, we must have that 
    \[\langle z_a,w\rangle=0 \quad \text{if}\quad |w|\neq 2\]
    and so it suffices to determine $\langle z_a,z_bz_c\rangle= -\eta_{a;bc}z_{a+b+c}$ for all $a,b,c\geq 1$.

    First note that
    \begin{align*}
        \langle z_a,z_1z_b\rangle &= (z_a,z_1z_b) - \langle z_1z_b,z_a\rangle\\
        &= -z_{a+b+1} -\langle z_1\succ z_b,z_a\rangle\\
        &= -z_{a+b+1}-\langle z_1,z_bz_a\rangle -\langle z_1,z_az_b\rangle.
    \end{align*}

    Next we find, using  the length-vanishing condition, we must have
    \begin{align*}
        0 &= \langle z_az_1z_bz_c,z_1\rangle\\
        &= \langle z_a\succ z_1z_bz_c,z_1\rangle - \langle \langle z_a,z_1z_b\rangle z_c, z_1\rangle\\
        &= \langle z_a,z_1\star z_1z_bz_c\rangle - \langle \langle z_a,z_1z_b\rangle z_c, z_1\rangle
        &= -\langle z_a, z_{b+2}z_c\rangle -\langle z_a,z_1z_{b+c+1}\rangle - \langle \langle z_a,z_1z_b\rangle z_c, z_1\rangle
    \end{align*}
    Then, using the previous observation, we find that
    \begin{equation*}\label{eq:gen-quasi-shuff-constraint}
    \langle z_a,z_{b+2}z_c\rangle = z_{a+b+c+1}+ \langle z_1,z_az_{b+c+1} + z_{b+c+1}z_a\rangle + \langle z_{a+b+1}z_c,z_1\rangle + \langle \langle z_1,z_az_b+z_bz_a\rangle z_c,z_1\rangle.
    \end{equation*}
    As such, for all $b\geq 3$, we can write $\langle z_a,z_bz_c\rangle$ in terms of $\langle z_1,z_mz_n\rangle$.

    Note that
    \[\langle z_a,z_bz_c\rangle = -z_{a+b+c}-\langle z_bz_c,z_a\rangle= -z_{a+b+c} -\langle z_b,z_cz_a\rangle -\langle z_b,z_az_c\rangle.\]
    If $b=1$, this is a linear combination of $\langle z_1,z_mz_n\rangle$. If $b=2$, we can use this identity again, along with the above reduction argument, to reduce to a linear combination of $\langle z_1,z_mz_n\rangle$ and $\langle z_2,z_2z_n\rangle$.

    As
    \[\langle z_2,z_2z_n\rangle =-z_{n+4} -\langle z_2,z_nz_2\rangle -\langle z_2,z_2z_n\rangle\]
    we can therefore reduce every bracket except $\langle z_2,z_2^2\rangle$ to a linear combination of brackets of the form $\langle z_1,z_mz_n\rangle$. 
    
    Noting that
    \[\langle z_2,z_2^2\rangle =-z_{6} -2\langle z_2,z_2^2\rangle \]
    we find that $\langle z_2,z_2^2\rangle = \frac{1}{3}(z_2,z_2^2)$. Thus, it remains to determine $\langle z_1,z_mz_n\rangle$.

    A similar calculation shows that
    \[\langle z_1,z_nz_1\rangle + \langle z_1,z_1z_n\rangle = \langle z_1\succ z_n,z_1\rangle = \langle z_1z_n,z_1\rangle = -z_{n+2}-\langle z_1,z_1z_n\rangle\]
    and hence
    \[\langle z_1,z_nz_1\rangle = -z_{n+2}-2\langle z_1,z_1z_n\rangle.\]
    Using this to simplify Equation (\ref{eq:gen-quasi-shuff-constraint}) with $a=1$, we find
    \[\langle z_1,z_{b+2}z_c\rangle = -\langle z_1,z_1z_{b+c+1}\rangle -\langle \langle z_1,z_1z_b\rangle z_c,z_1\rangle.\]
    Writing $\langle z_a,z_bz_c\rangle = -\eta_{a;b,c}z_{a+b+c}$, we therefore have
    \[-\eta_{1;b+2,c}z_{b+c+3} = \eta_{1;1,b+c+1}z_{b+c+3} + \eta_{1;1,b}\langle z_{b+2}z_c,z_1\rangle\]
    and so we can rearrange this to find
    \[\eta_{1;b+2,c} = \frac{\eta_{1;1,b}-\eta_{1;1,b+c+1}}{1+\eta_{1;1,b}}.\]
    Thus, it suffices to determine $\langle z_1, z_1z_n\rangle$ and $\langle z_1,z_2z_n\rangle$ for all $n$. Equation (\ref{eq:gen-quasi-shuff-constraint}) with $a=c=1$ allows us to determine $\eta_{1;1,n+2}$ in terms of $\eta_{1;1,n}$ for all $n$, and hence $\langle z_1, z_1z_n\rangle$ is entirely determined by 
    \[\langle z_1,z_1z_1\rangle = \frac{-z_3}{3}\quad  \text{and} \quad\langle z_1,z_1z_2\rangle=:-\alpha z_4\]
    Similarly, taking $a=1$, $c=2$, we find 
    $\langle z_1, z_1z_n\rangle$ is entirely determined by 
    \[\langle z_1,z_2z_1\rangle =-z_4-2\langle z_1,z_1,z2\rangle\quad  \text{and} \quad\langle z_1,z_2^2\rangle=:-\beta z_5\]

    Thus, any possible quasi-canonical structure is uniquely determined by the parameters $\alpha$ and $\beta$. Finally, we note that we must have
    \[\langle z_1z_2,z_1z_2z_1\rangle + \langle z_1z_2z_1,z_1z_2\rangle = z_7\]
    and
    \[\langle z_2z_1,z_2z_1z_2\rangle + \langle z_2z_1z_2,z_2z_1\rangle = z_8.\]
    Reducing both of these expressions in terms of $\alpha$ and $\beta$ leads to the equalities
    \begin{align*}
(10\beta-11)\alpha^2-22\alpha-10\beta+21 &=0\\
\alpha &=\frac{1}{2}
    \end{align*}
    and so $\alpha=\frac{1}{2}$, and $\beta=\frac{29}{30}$
    \end{proof}

In general, it may not be possible to determine all possible relations among the brackets $\langle z,w\rangle$ without some additional structure, which is why we cannot guarantee that this construction will always lead to genuine quasi-canonical structures. Similarly, it is not clear that a solution will exist at all. In this example, we have determined that at most one quasi-canonical structure is possible, but cannot verify that there are no additional clashing constraints on $\alpha$ and $\beta$. As such, a broader theory, more in line with the cohomological approach of the previous section, would be of much interest.

\vspace{1.5\baselineskip}
  \sloppy
  \bibliographystyle{abbrv}%
  \bibliography{bibliography.bib}%


\Needspace*{3\baselineskip}
\noindent
\rule{\textwidth}{0.15em}

{\noindent\small
School of Mathematics, Trinity College Dublin
Dublin 2, Ireland\
E-mail: \url{akeilthy@tcd.ie}\\
Homepage: \url{https://www.maths.tcd.ie/~keilthya}
}\vspace{.5\baselineskip}


\ifdraft{%
\listoftodos%
}

\end{document}
